\documentclass{article}

\usepackage[english]{babel}

\usepackage{amsfonts,amssymb,amsmath}

\usepackage{setspace}

\usepackage{centernot}

\def \erre {{\mathbb {R}}}
\def \gi {{\mathbb {G}}}

\def \elle {{\mathcal {L}}}

\def \acca{{\mathcal {H}}}
\def\j{\mathcal{J}}
\def \dive{\mathrm{div}}
 \def\sen{\mathrm{sen\,  }}

\def\erreu{{\erre^{ \scriptscriptstyle{N+1} }}}

\def\erren{{\erre^n}}

\def\inn{\ \mbox{ in }}
\def\ae{\mbox{\quad a.e.}}
\def\andd{ \quad\mbox{ and } \quad }
\def\c{\mathrm{const.}}
\def\forevery{\mbox{\quad for every }}
\def\ass{\quad \mbox{as}\quad}

\def\elleunoloc{{L^1_{\mathrm{loc}}}}

\def\cinfty {C^\infty(\erren,\erre)}
\def\ellep{L^p(\erren)}
\def\gradientea {{| \nabla_{\!\!\! A} u|}}

\newcommand{\tende}{\rightarrow}
\newcommand{\ttende}{\longrightarrow}
\newcommand{\enne} {\mathbb{N}}

\newcommand{\meno} {\backslash}

\newcommand{\inc} {\subseteq}
\newcommand{\inte} {\cap}

\newcommand\nequiv{\not\equiv}
\newcommand\eps{\varepsilon}
\newcommand{\frecciaf} {\longmapsto}

\newcommand{\vv}{{\overline V}} 
\newcommand{\pv}{{\partial V}}

\newtheorem{theorem}{Theorem}[section]
\newtheorem{proposition}[theorem]{Proposition}
\newtheorem{corollary}[theorem]{Corollary}
\newtheorem{lemma}[theorem]{Lemma}

\newtheorem{remark}[theorem]{Remark}

\numberwithin{equation}{section}

\def\proof{\noindent{\it Proof. }}
\def\endproof{\hfill $\Box$\par\vskip3mm}

\begin{document}

\title{${L^p}$-Liouville Theorems\\ for  Invariant Partial Differential Operators in ${\erre^n}$}

\author{{\sc{Alessia E. Kogoj and E. Lanconelli\thanks{Dipartimento di Matematica, Piazza di Porta San Donato 5, 40126 Bologna (Italy).
E-mail: alessia.kogoj@unibo.it, ermanno.lanconelli@unibo.it}}}
}
\date{ }
\maketitle

%\author{{\sc{Alessia E. Kogoj\thanks{Dipartimento di Matematica, Piazza di Porta San Donato 5, 40126 Bologna (Italy).}

%\begin{frontmatter}
%\author{Alessia E. Kogoj}
%\corref{cor1}}
%\thanks{alessia.kogoj@unibo.it}
%
%%\cortext[cor1]{Corresponding author.}
%
%\author{Ermanno Lanconelli}
%
%
%
%
%
%
%
%\address{Dipartimento di Matematica,\ Universit\`a di Bologna,\\
%Piazza di Porta San Donato, 5\  IT-40126 Bologna\  Italy}

\begin{abstract} 
We prove some $L^p$-Liouville theorems for hypoelliptic second order Partial Differential Operators
left translation invariant with respect to a Lie group composition law in
$\mathbb {R}^n$.  Results  for both solutions and subsolutions are given.
%Some $L^p$-Liouville theorems for several classes of Partial Differential Equations will be presented.
%The involved operators are left invariant with respect to Lie group composition laws in $\mathbb{R}^{n}$.
%Results for both solutions and sub-solutions will be given.
\end{abstract}

%\begin{keyword}Liouville Theorems \sep Invariant Partial Differential Operators  \sep  Hypoellliptic Operators on Lie groups 
%\MSC[2010]
%35B53   \sep   35R03 \sep   35H20  \sep 35H10 
%\end{keyword}

\section{Introduction}

In this paper we are concerned with $L^p$-Liouville properties for solutions and subsolutions to the equation
\begin{equation} \label{uno} 
\elle u = 0 \inn \erren,\end{equation}
where $\elle$ is a linear hypoelliptic second order Partial Differential Operator left translation invariant with respect to a Lie group in $\erren$. More precisely, the operator
$\elle$ in $\eqref{uno}$ is of the kind 
\begin{equation*} 
\elle: = \dive (A \nabla) + \langle b, \nabla \rangle,\end{equation*}
where $A=(a_{ij})_{i,j=1,\ldots,n}$ is a $n\times n$ symmetric matrix with real entries ${a_{ij}}$ in $C^\infty (\erren,\erre)$, 
$b=(b_1,\ldots, b_n)$ is a vector valued function with real components $b_j$ in  $C^\infty (\erren,\erre)$, and, as usual, $\dive$, $\nabla$ 
$\langle\ ,\  \rangle$ denote  Euclidean divergence, gradient and inner product in $\erren$.  We will assume, without further comments,
\begin{equation*} \langle A(x)\xi, \xi\rangle \geq 0 \quad \forall\ x, \xi \in \erren \andd \mathrm{trace\ }A (0) >0.\end{equation*} 
Our crucial assumptions on $\elle$ are the following ones.
\begin{itemize}
\item[(H1)] $\elle$ is hypoelliptic,  that is if $u$ is a distribution in a open set $\Omega\subseteq \erren$ and $\elle u$ is smooth in $\Omega$, then $u$ is smooth in $\Omega.$
\item[(H2)] There exists a Lie group $\gi=(\erren,\circ)$ such that $\elle$ is left translation invariant on $\gi$.
\end{itemize}  
For simplicity reasons we assume the {\it Lebesgue measure} in $\erren$ both  {\it left} and {\it right invariant} on $\gi$.
Throughout the paper we will  denote by $L^p$ the Lebesgue space $L^p(\erren,\erre)$.

 We recall that $\gi$ is said to be a {\it homogeneous Lie group} if the following property holds:  there exists a $n$-tuple of real numbers $\sigma=( \sigma_1,\ldots, \sigma_n)$, with $1 \le \sigma_1\le \ldots  \le \sigma_n$,
such that the dilation
\begin{eqnarray}\label{dilations} \delta_\lambda : \erre^n \longrightarrow \erre^n,\quad  \delta_\lambda ( x_1, \ldots, x_n) = (\lambda^{\sigma_1}x_1,\ldots,\lambda^{\sigma_n}x_n) \end{eqnarray}
is an  automorphisms of $\gi$, for every $\lambda >0.$  The real number 
$$Q= \sigma_1 + \cdots +  \sigma_n$$ is called the {\it homogeneous dimension} of $\gi$ w.r.t.  $(\delta_\lambda)_{\lambda >0}$.
 
If $\gi=(\erre^n,\circ, \delta_\lambda)$ is homogeneous then the Lebesgue measure in $\erren$ is right and left 
translation invariant on $\gi$ (see e.g. \cite{BLU}).

Aim of this paper is to prove the following theorems.
  \begin{theorem} \label{primo} Let $1\le p < \infty$ and let $u\in L^p$ be a smooth solution to 
  \begin{equation} \label{due} 
\elle u = 0 \inn \erren.\end{equation}
  Then $u\equiv 0.$
   \end{theorem}  
Nonnegative solutions to the equation  \eqref{due} satisfy also an $L^p$-Liouville property  for $0<p<1$. Indeed:

 \begin{theorem} \label{secondo} Let $0< p < 1$ and let  $u\geq 0$, $u^p \in L^1$, be a smooth solution to 
  \begin{equation*} 
\elle u = 0 \inn \erren.\end{equation*}
  Then $u\equiv 0.$
   \end{theorem}  
Theorem \ref{primo} extends to the subsolutions as follows.
   
   \begin{theorem} \label{terzo}  Let $u\in\elleunoloc $ be a solution to 
     \begin{equation*} 
\elle u \geq 0 \inn \erren, \mbox{in the weak sense of distributions.}\end{equation*}
 If  $u\in L^p$ for a suitable $p\in[1,\infty[,$ then $u\le 0$ a.e. in $\erren.$    \end{theorem}

When $\gi$ is a homogeneous group, Theorem \ref{terzo} takes the following sharp form.

\begin{theorem} \label{sharp}  Let  $\gi$ be a homogeneous Lie group with homogeneous dimension $Q\geq 3$. Assume 
$\elle$ is homogeneous of degree two w.r.t. the dilations in $\gi$.  Let $u\in\elleunoloc $ be a solution to 
\begin{equation*} 
\elle u \geq 0 \inn \erren, \mbox{in the weak sense of distributions.}\end{equation*}
If $u\in L^p$ for a suitable $p\in [ 1, 1 + {\frac{2}{Q-2}} ] $, then 
$$u\equiv 0 \ae\inn\erren.$$
Moreover, for every $p> 1 + {\frac{2}{Q-2}}$, there exists $u\in L^p, u\le 0$, $u\nequiv 0$, such that 
\begin{equation*} 
\elle u \geq 0 \inn \erren, \mbox{in the weak sense of distributions.}\end{equation*}
\end{theorem}  
Our proofs of the previous theorems are based on some devices that allow to obtain, as well, Liouville theorems 
for semilinear equations as the following one. We stress that this theorem does not requires the homogeneity  of $\gi$ and $\elle$.

\begin{theorem}\label{quinto}
Let $f:\erre\longrightarrow\erre$ be a $C^1$-increasing function such that \\
  \mbox{$f^{-1}(\{0\})=0.$}  Define  \begin{equation}\label{unoquattron} F:\erre\longrightarrow\erre,\quad F(t)= \int_0^t f(s)\ ds.\end{equation}
  Let $u\in C^2(\erren,\erre)$ be a classical solution to 
  \begin{equation}\label{unocinquen}\elle u = f(u) \inn \erren.\end{equation} 
  If $F(u)\in L^1(\erren)$ then $u\equiv 0.$
  \end{theorem}
 If in Theorem \ref{quinto} we choose  
  $f(t)=\lambda t$ or $f(t)=|t|^{p-1}t$, we obtain, respectively, the following corollaries.
 
  \begin{corollary} \label{sesto} Let $u\in C^2(\erren,\erre)\inte L^2(\erren)$ be such that  $$\elle u =\lambda u\inn\erren, \mbox{ with } \lambda\geq 0.$$ Then $u\equiv 0.$
  \end{corollary} 
  
    \begin{corollary}\label{settimo}
    Let  $1\le p <\infty$ and let $u\in L^{p+1}$ be a solution to       
    
    $$\elle u = |u|^{p-1} u  \inn \erren.$$
    Then $u\equiv 0.$

  \end{corollary}

\begin{remark}\label{noinfty}{\rm Theorem \ref{primo} does not hold, in general, if we assume $u\in L^\infty$ instead of $u\in L^p$ with $p<\infty$. Indeed, consider the Kolmogorov-type operator in $\erre^3=\erre_x^2 \times \erre_t$ 
$$\elle=\partial_{x_1}^2 + \left(x_1 - \frac{1}{2}x_2\right)  \partial_{x_1} +  \left( \frac{1}{2}x_1 - x_2\right)  \partial_{x_2} -\partial_t.$$ This operator satisfies (H1) and (H2), however, by a Priola and Zabczyk' s Theorem,  it has a bounded solution in $\erre^3$ which is not constant (see Remark \ref{priola} for details).}
\end{remark}
\begin{remark} {\rm When $\gi$ is a homogeneous group and $\elle$ is homogeneous w.r.t. the dilations of $\gi$, 
Theorem \ref{primo} and Theorem \ref{secondo} follow from a general Liouville-type theorem of Geller  \cite[Theorem 2]{geller_1983}. We want to stress that Geller' s Theorem also implies $L^\infty$-Liouville property for $\elle$ (if $\gi$ and $\elle$ are homogeneous). }
\end{remark}

\begin{remark} {\rm If the operator $\elle$ is homogeneous w.r.t. a group of dilations $(\delta_\lambda)_{\lambda>0}$ as in 
\eqref{dilations}, then Theorem \ref{primo} and Theorem \ref{secondo} hold only assuming hypothesis (H1). 

This follows from a result of Xuebo who extended Geller' s Theorem to homogeneous hypoelliptic
operators, not necessarily left invariant on a Lie group (see \cite[Theorem 1]{luo_xuebo_1997}).

}
\end{remark}

\begin{remark} {\rm We want to explicitly remark that Geller's and Xuebo's Theorems 
do not apply to subsolutions. }
\end{remark}

\begin{remark} {\rm  We say that $\elle$ satisfies the {\it one-side
Liouville property} if 
\begin{eqnarray*} \elle u = 0 \inn \erre^n , \ u\geq 0 \implies u=\mathrm{const.\  in\  } \erren.\end{eqnarray*} 
This property does not hold, in general, even for  left  translation invariant and homogenous operators. This is the case, e.g., of the classical heat operator in $\erren = \erre^{N+1}= \erren_x\times \erre_t$ 
$$\acca:= \Delta - \partial_t,\qquad  \Delta = \sum_{j=1}^n \partial_{x_j}^2,$$
which is hypoelliptic, invariant w.r.t. the euclidean translations and homogeneous of degree two  w.r.t.
the dilations 
\begin{eqnarray*} \delta_\lambda : \erreu \longrightarrow \erreu ,\quad  \delta_\lambda ( x,t) = (\lambda x, \lambda^2 t),\qquad  \lambda>0. \end{eqnarray*}
The function $u(x,t)=exp(x_1 + \cdots +  x_N + Nt)$ is a non-constant strictly positive solution to $\acca u=0$ in $\erre^{N+1}.$
Several classes of homogenous operators satisfying 
the one-side Liouville property have been presented in \cite{kogoj_lanconelli_2005, kogoj_lanconelli_2006, kogoj_lanconelli_2007}.}
\end{remark}

\begin{remark}{\rm $L^1$-Liouville Theorems for sub-Laplacians in suitable half spaces of stratified Lie groups
in $\erre^n$ have been proved by Uguzzoni \cite{uguzzoni_1999} and Kogoj \cite{kogoj_2014}.}
\end{remark}
\begin{remark}{\rm  When $\elle=\Delta$ is the classical Laplacian in $\erren$, Theorem \ref{quinto} is contained in \cite[Theorem 4.5]{caristi_dambrosio_mitidieri_liouv_2008}.}
\end{remark}

\begin{remark}{\rm For some kind of $L^p$-Liouville Theorems for sub-Laplacians on stratified Lie groups  we directly refer to the monograph \cite{BLU}, Chapter 5, Section 5.8.}
\end{remark}
\begin{remark}{\rm Liouville-type Theorems  based on suitable representation formulae for both solutions and subsolutions to some classes of higher order systems are contained in \cite{caristi_dambrosio_mitidieri_repr_2008}.}
\end{remark}

The remaining part of the present paper is organized as follows. In Section 2 we introduce a representation formula
which will play a crucial role in the proof of our Liouville-type Theorems. Some properties of the integral operators involved in the representation formula are proved in Section 3. Section 4 is devoted to the proofs of Theorems \ref{primo} and \ref{secondo}, while Theorem \ref{terzo} and \ref{sharp} are proved in Section 5 and 6, respectively. 
Section 7 contains the short proof of Theorem \ref{quinto}. Finally, in Section 8 we show some explicit examples of operators to which our results apply. 
\section{A representation formula}\label{2}
The assumptions  $A(0)\geq 0$ and $(\mathrm{tr\ } (A(0))>0$ imply the existence of a index $j\in\{1,\ldots,n\}$ such that 
$a_{jj}(0)>0$. For simplicity of notation we assume 
$$a_{11}(0)>0.$$
For $0<\eps<R$, let us define 

\begin{eqnarray*} V=V_{\eps, R} := D(R e_1, R+\eps) \inte  D(- R e_1, R+\eps),\end{eqnarray*}
where $e_1=(1,0,\ldots,0)$ and $D(\alpha, r)$ denotes the Euclidean ball with center $\alpha$ and radius $r$. If $R$ is sufficiently big and $\eps$ sufficiently smooth, $V$ satisfied the so called non-characteristic exterior ball condition at any point of its boundary. As a consequence, for every $x_0\in\partial V$ there exists a function $h(x_0,\cdot) \in C^2(\overline V, \erre)$ satisfying 

\begin{equation*}\begin{split} &h(x_0,\cdot)>0 \inn \overline{V} \meno {x_0} \andd h(x_0,x_0)=0,\\
&\elle_\eps h(x_0,\cdot)\le -1  \inn \overline{V} \mbox{\  for every\  }  \eps \in [0,1], \end{split}\end{equation*}
where $ \elle_\eps = \elle + \eps\Delta $ and $\Delta$ is the Laplace  operator in $\erren$.

The proof of this statement follows from very standard arguments (see e.g. \cite{BLU}, pages 383, 384 and 387).

The existence of barrier functions  $h(x_0,\cdot)$ implies the following Picone-type estimate: there exists a constant $C>0$, independent of $\eps\in[0,1],$ such that 
\begin{eqnarray} \label{piestimate} \sup_{\overline V} |u|\le \sup_{\partial V} |u| + C \sup_{\overline V} |\elle_\eps u|,\qquad \forall\ u \in C^2(\overline V, \erre)\end{eqnarray} and for every $\eps\in[0,1].$

We also have the following Picone Maximum Principle:\\ if $u\in C^2(V, \erre)\inte C(\overline V, \erre)$ satisfies
$$\elle_\eps u\geq 0 \inn V\andd u|_{\partial V} \le 0,$$
then $u\le 0$ in $V$ (see, e.g. \cite{lanconelli_2010}).

The hypoellipticity  of $\elle$, the estimate \eqref{piestimate} and the boundary barrier functions $h(x_0,\cdot)$ allow to prove the solvability of the boundary value problem 
\begin{equation} \label{bvp}
\begin{cases}
\elle u= -f  \inn V   \\ 
u|_{\partial U}  = \varphi 
\end{cases}
\end{equation}
with an {\it elliptic regularization procedure}. 
\begin{proposition} For every $f\in C^{\infty}(\overline{V}, \erre)$ and for every $\varphi \in C(\partial V, \erre)$ the boundary value problem \eqref{bvp} has a unique solution $u\in C^\infty(V,\erre)\inte C(\overline V, \erre).$ This solution satisfies the estimate 
\begin{eqnarray*}  \sup_{ V} |u|\le \sup_{\partial V} |\varphi| + C \sup_{ V} |f|,\end{eqnarray*} 
where $C>0$ does not depend on $u, \varphi$ and $f$. Moreover, if $f\geq 0$ and $\varphi\geq 0$, then $u\geq 0.$ 
\end{proposition}

\proof It follows, along standard lines, as in \cite[Theorem 5.2] {bony_1969} and \cite[pages 383--387 ]{BLU}.
\endproof
We denote by $G(f)$ the solution of \eqref{bvp} with  $\varphi = 0$ and by $H(\varphi)$ the solution of \eqref{bvp}
with $f=0$. Then, letting $C_0(\overline V, \erre)=\{ u\in C(\overline V, \erre)\ |\ u=0 \mbox{ on } \partial V\},$
the operators
\begin{equation*}\begin{split} & G: C^\infty(\overline V, \erre)\ttende C_0(\overline V, \erre)\\ \andd \\
& H: C(\partial V, \erre)\ttende C^\infty(V,\erre)\inte C(\vv, \erre)
\end{split} 
\end{equation*}
are {\it linear, nonnegative} and satisfy 

\begin{equation*}\begin{split} & \sup_\vv |G(f)|\le C\sup_\vv |f| \qquad \forall f\in C^\infty(\vv,\erre),\\
&  \sup_\vv |H(\varphi)|\le \sup_{\partial V}  |\varphi|.
\end{split} 
\end{equation*}
Then, $G$ can be continued to a linear, nonnegative and continuous operator, still denoted by $G$,
$$G: C(\overline V, \erre)\ttende C_0(\overline V, \erre).$$
Let us now consider the functionals 
\begin{equation*}\begin{split} & C(\vv,\erre) \ni f \frecciaf G(f)(0) \in\erre \\  \andd \\
& C(\partial V ,\erre) \ni \varphi  \frecciaf H(\varphi)(0) \in\erre.
\end{split} 
\end{equation*}
They are linear, nonnegative and continuous. Then there exist two nonnegative Radon measures $\nu$ and $\mu$, respectively on $\vv$ and $\pv$ such that 
\begin{equation*}\begin{split} &  G(f)(0) = \int_\vv f \ d\nu \qquad \mbox{ for every } f\in C(\vv, \erre) \\  \andd \\
&  H(\varphi)(0) = \int_\pv \varphi \ d\mu \qquad \mbox{for every } \varphi \in C(\pv, \erre).
\end{split} 
\end{equation*}
Then, the following proposition holds
\begin{proposition}\label{propformularapp} For every $u\in C^2(\vv,\erre)$ we have 
\begin{equation} \label{formularapp} u(0)= \int_\pv u \ d\mu - \int_\vv \elle u \ d\nu.
\end{equation}
\end{proposition}
\proof  We first assume $u\in C^\infty(\vv,\erre)$. Let us put $f=\elle u$. By the very definition of $G$, the function $v:=G(f)$ satisfies: $v\in C^\infty (V,\erre) \inte C(\vv,\erre)$ and 
\begin{equation*} 
\begin{cases}
\elle u= -f = - \elle u   \inn V,  \\ 
v|_{\partial V}  = 0.
\end{cases}
\end{equation*}
Then, $\elle (v+u)=0$ and $(v+u)|_\pv = u|_\pv,$ so that 
$$v+u= H(u|\pv),$$
i.e., 
$$u= H(u|\pv) - G(\elle u).$$
Hence 
\begin{equation*}\begin{split}
u(0) &= H(u|\pv)(0) - G(\elle u)(0) \\
&= \int_\pv u d \mu - \int_\vv \elle u \ d\nu.\end{split}
\end{equation*}
Then \eqref{formularapp} holds true if $u\in C^\infty(\vv, \erre)$. On the other hand, if $u\in C^2(\vv, \erre)$, there exists
a sequence $(u_n)$, with $u_n\in C^\infty(\vv, \erre)$, such that 

\begin{equation*}\begin{split}
& u_n\ttende u \mbox{ \quad uniformly  on } \pv \\
\andd \\
&\elle u_n\ttende \elle u \mbox{ \quad uniformly  on } \vv .\end{split}
\end{equation*}
On the other hand, for what we have already proved, 
$$u_n(0) = \int_\pv  u_n d \mu - \int_\vv \elle u_n\ d\nu\qquad \forall\ n\in \enne.$$
Letting $n$ go to infinity we obtain \eqref{formularapp}.
\endproof
We are now ready to state and prove the main result of this section.

\begin{theorem} \label{rappresentazione} For every $v \in C^2(\erren, \erre)$ we have 
\begin{eqnarray*} v(x)= M(v)(x) - N (\elle v) (x)\qquad \forall \ x\in\erren,\end{eqnarray*}
where,
\begin{eqnarray}\label{M} M(u)(x)= \int_\pv u(x \circ y)\ d \mu(y),\end{eqnarray}
and 
\begin{eqnarray} \label{N}N(f)(x)= \int_\vv f(x \circ y)\ d \nu(y).\end{eqnarray}
Here $\circ$ denotes the composition law of $\gi.$
\end{theorem}
\proof  Let $x\in\erren$ be fixed and consider the function 
$$u: \erren \ttende \erre,\qquad u(y)= v(x\circ y).$$
Obviously $u\in C^2(\erre^n,\erre)$ and, since $\elle$ is left translation  invariant,
$$\elle u (y) = (\elle v) (x \circ y) \mbox{ for every } y\in\erren.$$ Then, by Proposition \ref{propformularapp},
\begin{equation*} \begin{split} v(x)= u(0)=& \int_\pv u(y)\ d\mu(y) - \int_\vv \elle u(y)\ d\nu(y)\\
&= \int_\pv v(x \circ y)\ d\mu(y) - \int_\vv (\elle v) (x \circ y)\ d\nu(y)\\
&= M(v)(x) - N (\elle( v)) (x).
 \end{split}
\end{equation*}
\endproof
\section{Some properties of the operators $M$ and $N$}
In this Section we prove some properties of the operators $M$ and $N$ defined in \eqref{M} and \eqref{N}, respectively.

We start with the following lemma.
\begin{lemma}\label{elleunoelleuno} Let $v$ be a continuous and $L^1$-function in $\erren.$\\ Then $M(v) \in L^1(\erren)$ and 
\begin{eqnarray} \label{uffi} \int_\erren v(x)\ dx = \int_\erren M(v)(x)\ dx.
\end{eqnarray}
\end{lemma}
\proof 
It follows from Fubini Theorem and the invariance of the Lebesgue measure on $\gi$. Indeed:
\begin{equation}\label{chissa}\begin{split}\int_\erren M(v)(x)\ dx &= \int_\erren \left( \int_\pv v(x\circ y)\ d\mu(y)\right)\ dx\\
&= \int_\pv \left( \int_\erren v(x\circ y)\ dx \right)\  d\mu(y)\\
&= \left(\int_\erren v(x)\ dx \right) \left(\int_\pv  \  d\mu(y)\right).
\end{split}\end{equation} 
On the other hand, by Proposition  \ref{propformularapp} applied to the function $u\equiv 1$, we have 
$$1 = \int_\pv d\mu.$$
Using this information in \eqref{chissa}, we obtain \eqref{uffi}.
\endproof
Regarding the operator $N$ we have:
\begin{lemma}\label{unduetre}Let $f\in C(\erren, \erre).$ Then the following statements hold.
\begin{itemize}
\item[ (i)] $N(f) \geq 0$ if $f\geq 0$;

\item[(ii)] if $f\geq 0$ and $N(f)\equiv 0$, then $f\equiv 0$;

\item[(iii)] $N(f)\in C(\erren,\erre)$.

\end{itemize}
\end{lemma} 
\proof\\ {\it (i)} It is obvious.
\\ {\it (ii)} Let $f\geq 0$ and $N(f)\equiv 0$ in $\erren$. Assume, by contradiction, $f\nequiv 0$. Then, there exists 
$x_0 \in\erren$ such that $f(x_0)>0.$ Since $f$ is continuous, there exists an open set $\Omega\ni x_0$ such that 
$f(x)>0$ for every $x\in\Omega$. It follows 
$$f(z\circ y)>0\qquad \forall\ z, y \in \erren \ : \ z\circ y\in\Omega.$$
As a consequence, for every $z\in\erren$ we have
$$0 = N(f)(z)= \int_\vv f(z\circ y)\ d\nu(y) \geq \int_{\vv\inte(z^{-1}\circ\Omega)} f(z\circ y)\ d\nu(y),$$
so that, since $f(z\circ y)>0$ for every $y\in z^{-1}\circ \Omega,$ we get 
\begin{eqnarray} \label{paperinik} \nu(\vv\inte (z^{-1} \circ \Omega))=0\qquad \forall z\in\erren.
\end{eqnarray} Since 
$$\bigcup_{z\in\erren} (z^{-1} \circ \Omega) = \erren,$$ 
from \eqref{paperinik} we obtain 
$$\nu(\vv)=0.$$
As a consequence, by Proposition \ref{propformularapp},
$$u(0)=\int_\pv u \ d\mu\qquad \forall u\in C^2(\vv,\erre).$$
In particular:
$$u(0)=0\qquad \forall u\in C_0^\infty(V,\erre),$$ which is absurd. This completes the proof of {\it (ii)}.
\\ {\it (iii)} Since $f$ is continuous and $\vv$ is compact, for every $z_0\in\erren$ we have:
$$\sup_{y\in\vv} |f(z\circ y) - f(z_0 \circ y)| \ttende 0 \mbox{ as } z\ttende z_0.$$
Then
\begin{equation*}\begin{split} N(f)(z_0) = & \int_\vv f(z_0 \circ y)\ d\nu(y) = \int_\vv \lim_{z\tende z_0} f(z\circ y)\ d\nu(y)\\
&=  \lim_{z\tende z_0} \int_\vv  f(z\circ y)\ d\nu(y)\\
&=  \lim_{z\tende z_0} N(f)(z).
\end{split}\end{equation*} This proves the continuity of $N(f).$ 

\endproof 
\section{Proof of Theorems \ref{primo} and \ref{secondo}}\label{dimprimoesecondo}
We start with the following elementary lemma.
\begin{lemma} \label{elementarylemma} Let $F:\erre\ttende\erre$ be a $C^2$-function and let $u\in C^2(\Omega, \erre), \Omega \subseteq \erren,$ open. Then $$v:=F(u)$$ is a real $C^2$-function in $\Omega$ such that
\begin{equation*}\elle v = F' (u)\elle u + F''(u) | \nabla_{\!\!\! A} u |^2, \end{equation*} where

 \begin{equation*} | \nabla_{\!\!\! A} u |^2:= \langle A \nabla u, \nabla u\rangle.
\end{equation*}
\end{lemma} 
\proof We show the elementary computations for reader convenience. We have:
\begin{equation*}\begin{split} \elle v = & \dive (A \nabla (F(u)) + \langle b, \nabla (F (u))\rangle \\
= & \dive (F'(u) A \nabla u)) + \langle b, \nabla  u\rangle F'(u)  \\
= & F'(u) ( \dive (A \nabla u) + \langle b, \nabla  u\rangle) +  F''(u) \langle A \nabla u, \nabla u\rangle \\
= & F'(u)\elle u + F''(u)|\nabla_{\!\!\! A} u|^2.
\end{split}
\end{equation*}\endproof
To prove our theorems we need another lemma.
\begin{lemma} \label{anotherlemma}  Let $\Omega \inc \erren$ be open and connected and let $v\in C^2(\Omega, \erre)$ be such that 
\begin{equation}\label{propagazione} |\nabla_{\!\!\! A} u|^2 = 0 \andd \langle b, \nabla u \rangle = 0 \inn \Omega.\end{equation} 
Then $u=\mathrm{const.}$ in $\Omega$.\end{lemma} 
\proof Let us denote by $X_1, \ldots, X_n$ the vector fields constructed with the columns of the matrix $A$, i.e., 
\begin{equation*} X_j=\sum_{k=1}^{n} a_{k,j} \partial_{x_k},\qquad j=1,\ldots,n.\end{equation*}
Let us also put
\begin{equation*} X_0=\sum_{k=1}^{n} b_{k} \partial_{x_k}.\end{equation*}
Then, assumption \eqref{propagazione} can be written as follows 
\begin{equation*} X_j u = 0 \inn \Omega \mbox{\quad for every } j=1,\ldots,n.\end{equation*}
As a consequence, 
\begin{equation} \label{trestarn}Y u = 0 \inn \Omega \qquad\forall \ Y\in \mathrm{Lie}  \{X_0, X_1, \ldots, X_n\}. \end{equation}
On the other hand, since $\elle$ is hypoelliptic, 
\begin{equation*} \mathrm{rank\ } \mathrm{Lie}  \{X_0, X_1, \ldots, X_n\}(x)=N \qquad\forall\ x\in \Omega^*, \end{equation*}
where $\Omega^*$ is an open subset of $\Omega$ such that $\overline{\Omega^*}\inc\Omega.$ Then, for every 
$x\in\Omega^*$ and for every $i\in\{1,\ldots,n\}$ there exist $Y_1,\ldots,Y_n \in \mathrm{Lie} \{X_0,X_1,\ldots, X_n\}$ and real constants $c_i^{(1)} (x), \ldots, c_i^{(n)} (x)$ such that 
$$\partial_{x_i}= \sum_{j=1}^n c_i^{(j)} Y_j.$$
Thus, from \eqref{trestarn}, we obtain  
$$\partial_{x_i} u(x)=0\qquad\forall\ x\in\Omega^*,\ \forall\ i =1,\ldots,n.$$
Since $u\in C^1(\Omega, \erre)$ and $\Omega^*$ is dense in $\Omega$, this implies 
$$\nabla u\equiv 0\inn \Omega.$$
Then $u$ is constant in $\Omega$.
\endproof 
A key tool in the proof of our Liouville theorems is given by the following proposition.
\begin{proposition} \label{0}Let $v\in C^2(\erren, \erre).$ If 
$$ v\in L^1(\erren) \andd \elle v\geq 0 \inn\erren$$
then
$$\elle v=0\inn\erren.$$
\end{proposition}
\proof The representation formula of Theorem \ref{rappresentazione} gives 
\begin{eqnarray*} v= M(v) - N (\elle v).\end{eqnarray*}
Now, being $ v\in L^1(\erren)$, Lemma \ref{elleunoelleuno} implies
$ M(v)\in L^1(\erren)$ and \begin{eqnarray*}  \int_\erren v(x)\ dx = \int_\erren M(v)(x)\ dx.
 \end{eqnarray*}
Then, $N(\elle v) \in L^1(\erren)$ and
 \begin{eqnarray*}  \int_\erren N (\elle v)(x)\ dx =0.
 \end{eqnarray*}
 Since $\elle v \geq 0$, by Lemma  \ref{unduetre}-{\it (i)}, $N(\elle v) \geq 0$ and the last integral identity implies 
 $N(\elle v) = 0$ a.e. in $\erren$. 
 On the other hand, $\elle v$ is continuous and, by Lemma \ref{unduetre}-{\it (iii)}, $N(\elle v)$ is continuous. Therefore
 $$N(\elle v)=0\inn\erren,$$so that, by Lemma \ref{unduetre}-{\it (ii)},
 $$\elle v =0\inn\erren.$$
\endproof

We are now ready to prove Theorem \ref{primo}.

{\noindent{\it Proof of Theorem \ref{primo}. }} 
Let $u\in L^p(\erren)$, $1\le p <\infty,$ be a solution to $\elle u=0$ in $\erren$ and assume, by contradiction, $u\nequiv 0$. Define
   $$v:= F(u)$$ with $$F:\erre\longrightarrow\erre,\quad  F(t)=(\sqrt{1+t^2} - 1)^p .$$
  Elementary computations show that $F\in C^2(\erre,\erre),$
  
 \begin{equation}\label{tretre} 0\le F(t)= \left(\frac{t^2}{ \sqrt{1+t^2} +1}\right)^p\le |t|^p\end{equation} 
 and
  \begin{equation}\label{trequattro} F''(t)>0\quad\forall t \neq 0.\end{equation}
Since $\elle$ is hypoelliptic the function $u$ is smooth and, by Lemma \ref{elementarylemma},
we have
\begin{equation}\label{trecinque}\elle(F(u)) = F' (u)\elle u + F''(u) | \nabla_{\!\!\! A} u |^2 =  F''(u) | \nabla_{\!\!\! A} u |^2.\end{equation} 
On the other hand $F(u)\in L^1(\erren)$, since, by \eqref{tretre},
$$ 0\le F(u)\le |u|^p\andd u\in L^p(\erren).$$
From Proposition \ref{0} it follows 
$$\elle(F(u))=0\inn\erren,$$
so that, keeping in mind \eqref{trecinque}, 
$$F''(u) | \nabla_{\!\!\! A} u |^2 =0.$$
Then, by \eqref{trequattro}, 
\begin{equation}\label{tresei} | \nabla_{\!\!\! A} u |^2=0\inn \Omega_0:=\{x\in\erren\ |\ u(x)\neq 0\}.\end{equation}
$\Omega_0$ is an open subset of $\erren$ which is nonempty because we are assuming $u\nequiv 0$. 
Since $A\geq 0$ from \eqref{tresei} we obtain 
$$A\nabla u=0\inn \Omega_0,$$
so that
$$\dive(A\nabla u) =0\inn \Omega_0.$$
As a consequence, keeping in mind that $\elle u=0$, 
\begin{equation}\label{tresette} \langle b, \nabla u\rangle  =0\inn \Omega_0.
\end{equation}
Identities \eqref{tresei} and \eqref{tresette} and Lemma \ref{anotherlemma}  imply
\begin{equation*} u=\c \mbox{ on every connected component of } \Omega_0.
\end{equation*}
Let $O$ be one of the connected component of $\Omega_0$.  If $O=\erren$ we have $u=\c$ in $\erren$, so that, since $u\in L^p(\erren)$, $u=0$ in $O$.  If $O\neq\erren$,  then $\partial O\neq \emptyset$ and $u=0$ on $\partial O$. Being $u=\c$ in $O$, this implies $u=0$ in $O$.

Thus, we have proved that $u=0$ on every connected component of $\Omega_0$, that is
$$u=0\inn \Omega_0,$$
in contradiction with the definition of $\Omega_0.$
\endproof
The previous argument  can be easily adapted to prove Theorem \ref{secondo}. 
\\{\noindent{\it Proof of Theorem \ref{secondo}. }} 
Let $u\geq 0$ be a solution to $\elle u=0$ such that $u^p\in L^1(\erren)$ for a suitable $p\in ]0,1[.$ Define $$v:= F(u)$$with 
$$F: [0, \infty[ \longrightarrow\erre,\quad  F(t)=(1 + t)^p - 1.$$
The function $F$ is smooth, 
  
 \begin{equation*} 0\le F(t)\le t^p \forevery t\geq 0\end{equation*} 
 and
  \begin{equation*} F''(t)< 0 \forevery t\geq 0.\end{equation*} 
Since $\elle$ is hypoelliptic the function $u$ is smooth. By Lemma \ref{elementarylemma},
we have
\begin{equation*}\elle(F(u)) = F' (u)\elle u + F''(u) | \nabla_{\!\!\! A} u |^2 =  F''(u) | \nabla_{\!\!\! A} u |^2.\end{equation*} 
Hence $\elle(F(u))\le 0.$ On the other hand $F(u)\in L^1(\erren)$, since
$$ 0\le F(u)\le u^p\andd u^p \in L^1(\erren).$$
Then by Proposition \ref{0}, 
$$\elle(F(u))=0\inn\erren,$$
so that, 
$$F''(u) | \nabla_{\!\!\! A} u |^2 =0 \inn \erren.$$
Being $F''(u(x))<0$ at any point, from this last identity we obtain 

\begin{equation*}| \nabla_{\!\!\! A} u |^2=0\inn \erren.\end{equation*}
Now, arguing as in the proof of Theorem \ref{primo}, we obtain 
$$u=\c\inn\erren,$$ so that, since $u^p\in(\erren),$ 
$$u=0 \inn\erren.$$
\section{Proof of Theorem \ref{terzo}}
We start by proving Theorem \ref{terzo} in the case of $u$ smooth. Thus, let \mbox{$u\in C^\infty(\erren,\erre)$} be such that 
$$\elle u \geq 0\inn\erren\andd u\in L^p(\erren),\  1\le p < \infty.$$
We want to prove that $u\le 0$ in $\erren.$ Arguing by contradiction, we assume
$$\Omega_0:=\{ x\in\erren\ | \  u(x)>0\}\neq\emptyset.$$
Let us consider the function 
\begin{equation*} F:\erre\ttende\erre,\quad F(t) = \begin{cases} 0 & \mbox{\  if } t\le 0, \\ 
(({1 + t^4})^{\frac{1}{4}} - 1)^{p}  & \mbox{\  if } t> 0.
\end{cases}\end{equation*}
Then:
\begin{itemize}
\item[$(i)$] $F\in C^2(\erre,\erre);$
\item[$(ii)$]  $F$ is increasing and convex;
\item[$(iii)$] $F''(t)>0$ if $t>0$;
\item[$(iv)$] $0\le F(t)\le t^p$ for every $t>0.$
\end{itemize}
We let $$v:=F(u).$$
From the properties of $F$ we get that $v\in C^2(\erre,\erre)$ and $0\le v \le C |u|^p$, so that $v\in L^1(\erren).$
Moreover 
$$\elle v= F'(u) \elle u + F''(u)\gradientea^2 \geq 0.$$
Then, by Proposition \ref{0}, $\elle v = 0$ hence, in particular 
\begin{equation}\label{uguale0} F''(u)\gradientea^2 = 0 \inn \erren.\end{equation}
Since $F''(u(x))>0$ for every $x\in\Omega_0,$ \eqref {uguale0} implies 
\begin{equation*} \gradientea^2 = 0 \inn \Omega_0.\end{equation*}
Starting from this identity and arguing as in the proof of Theorem \ref{primo} (see Section \ref{dimprimoesecondo}), we obtain 
\begin{equation*} u = 0 \inn \Omega_0, \end{equation*}
in contradiction with the definition of $ \Omega_0$. This proves Theorem \ref{terzo} in the case $u$ smooth. We will remove this restriction by using the following lemma.
\begin{lemma} Let $u\in\ellep, 1\le p < \infty,$ be such that $\elle u \geq 0$ in $\erren$ in the weak sense of distributions. Then there exists a sequence of functions $(u_k)$ such that 
\begin{itemize}
\item[$(i)$] $u_k \in \cinfty$ for every $k\in\enne$;
\item[$(ii)$]  $\elle u_k\geq 0$ in $\erren$  for every $k\in\enne$;
\item[$(iii)$]  $\elle u_k\in\ellep$ in $\erren$  for every $k\in\enne$;
\item[$(iv)$] $u_k\ttende u$ in $\elleunoloc(\erren)$.
\end{itemize}
\end{lemma}
\proof The proof is quite standard. We give it in the details for reader convenience. Let $\eps >0$ be fixed and choose 
a function $\mathcal{J}_\eps \in C_0^\infty(\erren,\erre)$ such that  $\mathrm{supp}\ \j_\eps\inc D(0,\eps),$ 
$\int_\erren \j_\eps (y)\ dy=1$,  $ \j_\eps\geq 0$. Define 
$$\hat u_\eps:\erren\ttende\erre, \hat u_\eps(x)=\int_\erren u(y\circ x) \j_\eps(y)\ dy.$$
A change of variable in the integral gives 
$$\hat u_\eps(x)=\int_\erren u(z) \j_\eps(z\circ x^{-1})\ dz,$$
showing that 
$u_\eps \in C_0^\infty(\erren,\erre).$ Moreover, for every $\varphi\in C_0^\infty(\erren), \varphi\geq 0$, we have 
($\elle^*$ = formal adjoint of $\elle$)
\begin{equation*}\begin{split} \int_\erren \hat u_\eps (x) \elle^* \varphi (x)\ dx &=\int_\erren\left(\j_\eps(y) \int_\erren  u(y\circ x) \elle^*\varphi(x)\ dx\right)\ dy \\ &=\int_\erren \j_\eps(y) \left(\int_\erren u(z) (\elle^*\varphi) (y^{-1}\circ z) \ dz\right)\ dy  \\
&=  \mbox{ (since $\elle^*$ is left translation invariant)}
 \\ &\phantom{=} \int_\erren \j_\eps(y) \left(\int_\erren  u(z) (\elle^*(\varphi (y^{-1}\circ z)) \ dz\right)\ dy.
 \end{split}
\end{equation*}
Since $\elle u\geq 0$ in the weak sense of distributions the inner integral at the last right hand side is $\geq 0$. Therefore
\begin{equation*} \int_\erren \hat u_\eps (x) \elle^* \varphi (x)\ dx \geq 0 \qquad\forall\ \varphi\in C_0^\infty(\erren),\ \varphi\geq 0.\end{equation*} Since $\hat u_\eps$ is smooth we can integrate by parts at the left hand side, getting 
\begin{equation*} \int_\erren \elle \hat u_\eps (x)  \varphi (x)\ dx \geq 0 \qquad\forall\ \varphi\in C_0^\infty(\erren),\ \varphi\geq 0.\end{equation*}
Thus 
\begin{equation*} \elle \hat u_\eps \geq 0 \inn \erren.\end{equation*}
Moreover 
\begin{equation*}\begin{split}\int_\erren |\hat u_\eps (x)|^p\ dx &  \le \int_\erren  \left( \int_\erren |\hat u (y\circ x)|^p \j_\eps(y)\ dy \right) \ dx \\ &
= \int_\erren  \j_\eps(y)  \left( \int_\erren |\hat u (y\circ x)|^p \ dx \right) \ dy
\\ &
=  \int_\erren |\hat u ( x)|^p \ dx.
\end{split}\end{equation*}
Hence
\begin{eqnarray*} \hat u_\eps \in L^p(\erren).
\end{eqnarray*}
Finally, for every fixed compact set $K\inc\erren$, 

\begin{equation*}\begin{split}\int_K |\hat u_\eps (x) - u(x)|\  dx &  \le \int_\erren \j_\eps(y)  \left( \int_K |u (y\circ x) - u(x)|\ dx \right) \ dy\\
&  \le \sup_{y\in D(0,\eps)} \int_K |u (y\circ x) - u(x)|\ dx 
\\ & := \omega_K(u,\eps).
\end{split}\end{equation*}
On the other hand, being $u\in\elleunoloc(\erren)$,
$$\omega_K(u,\eps)\ttende 0 \quad \mbox{as}\quad {\eps\ttende 0}.$$
Therefore, a sequence $(u_k)_{k\in\enne}$ satisfying $(i)-(iv)$ can be constructed by choosing $u_k=\hat u_{\frac{1}{k}}.$ \endproof
We are ready to complete the proof of Theorem \ref{terzo}. 

Let  $u\in\ellep, 1\le p<\infty$, be such that $\elle u\geq 0$ in the weak sense of distributions. By the previous Lemma there exists a sequence $(u_k)_{k\in\enne}$ of smooth functions such that $\elle u_k\geq 0, u_k\in\ellep$ and $u_k\ttende u$ as $k\ttende\infty$ in $\elleunoloc(\erren).$  For what proved in the first part of this section, $u_k\le 0$ in $\erren$ for every $k\in\enne.$ This implies $u\le 0$ a.e. in $\erren$, and completes the proof.

\section{Proof of the Theorem \ref{sharp}}
We need several prerequisites. First of all, the assumptions on $\elle$ and $\gi$ imply the existence of a fundamental solution $$\Gamma:\erren\ttende [0,\infty]$$ such that 
\begin{itemize}
\item[$(i)$] $\Gamma\in \elleunoloc(\erren), \Gamma\in C^\infty(\erren \meno \{0\})$ and
$$\Gamma(x)\ttende 0 \ass x\ttende\infty;$$
\item[$(ii)$]  $\int_\erren \Gamma(x) \elle^* \varphi (x)\ dx = - \varphi(0)\qquad\forall \varphi\in C_0^\infty(\erren);$
\item[$(iii)$]  $\Gamma(\delta_\lambda(x)) =\lambda^{2-Q}\ \Gamma(x)\qquad\forall\ x\in \erren\meno \{0\},\ \forall\ \lambda>0$
\end{itemize}
(see e.g. \cite{folland_1975}).\\
Given $f\in C_0^\infty(\erren,\erre)$ we let 
$$\Gamma \ast f(x):=\int_\erren \Gamma(y^{-1}\circ x) f(y)\ dy = \int_\erren \Gamma(z) f(x\circ z^{-1})\ dz.$$
From $(i)$ it follows:
$$\Gamma\ast f \in C^\infty(\erren,\erre)\andd \Gamma\ast f(x)\ttende 0 \ass x\ttende\infty.$$
Moreover, as an elementary computation shows,
$$\int_\erren(\Gamma\ast f)(x) \elle^* \varphi(x)\ dx = - \int_\erren f(y)\varphi(y)\ dy\qquad\forall\ \varphi\in C_0^\infty(\erren,\erre).$$ Hence, 
$$\elle(\Gamma\ast f)=- f.$$
The operator $\elle$ satisfies the following Maximum Principle on $\erren$.
\begin{proposition} \label{maxprinc61} Let $u\in C^2(\erren,\erre)$ be such that 
\begin{equation}\label{seidue}\elle u\geq 0 \inn\erren\andd \limsup_{x\ttende\infty} u(x) \le 0.\end{equation}
Then $u\le 0$ in $\erren.$ 
\end{proposition} 
\proof  We have already remarked the existence of a bounded neighborhood $V$ of the origin on which $\elle$ satisfies the Picone Maximum Principle: \begin{center}if $v\in C^2(\erren,\erre)\inte C(\overline V, \erre), \elle v\geq 0$ in $V$ and $v|_{\partial V} \le 0,$ then $v\le0$ in $V$\end{center}
(see Section \ref{2}). Let $\eps >0$ be arbitrarily fixed and define 
$$v_\lambda(x):=u(\delta_\lambda(x))-\eps,\quad x\in\erren,\ \lambda>0.$$
The second assumption in \eqref{seidue} implies the existence of $\lambda_\eps>0$ such that 
$$v_\lambda(x)\le 0\qquad\forall\ x\in\partial V,\ \forall\ \lambda>\lambda_\eps.$$
Moreover, since $\elle$ is $\delta_\lambda$-homogeneous of degree two:
$$\elle v_\lambda(x)=\lambda^2(\elle u)(\delta_\lambda(x))\geq 0 \qquad\forall\ x\in\erren,\ \forall \lambda >0.$$ As a consequence, by Picone Maximum Principle on $V$,
$$v_\lambda\le 0\inn V\quad\forall \lambda>\lambda_\epsilon,$$which means
\begin{equation}\label{seitre} u(\delta_\lambda(x))\le\eps\qquad\forall\ x\in V,\ \forall \lambda>\lambda_\eps.\end{equation}
On the other hand,  since $V$ is a neighborhood of the origin 
$$\bigcup_{\lambda>\lambda_\eps} \delta_\lambda(V)=\erren.$$
Together with \eqref{seitre} this implies 
$$u\le\eps\inn\erren,\qquad\forall\ \eps>0.$$
Hence $u\le 0$ in $\erren$.
\endproof 
As an application of the previous proposition, we prove the positivity of $\Gamma$.
\begin{corollary}
It is
\begin{equation}\label{seiquattro} \Gamma(x)\geq 0\qquad\forall x\in\erren\meno\{0\}.\end{equation}
\end{corollary}
\proof For every $f\in C_0^\infty(\erren,\erre), f\le 0,$ we have 
$$\elle(\Gamma\ast f)=-f\geq 0\andd \Gamma\ast f|_{\infty}=0.$$
Then, by the previous theorem, $\Gamma\ast f \le 0$ in $\erren.$ In particular 
$$\Gamma\ast f(0)=\int_\erren \Gamma(y^{-1}) f(y)\ dy \le 0\qquad\forall\ f\in C_0^\infty(\erren,\erre),\ f\le 0,$$
from which \eqref{seiquattro} follows, since $\Gamma$ is smooth out of the origin.
\endproof
\noindent {\bf Note.}  If we agree to let $$\Gamma(0):=\liminf_{x\ttende 0} \Gamma(x),$$
then $\Gamma:\erren\ttende [0,\infty]$ is {\it lower semicontinuous}.\\
Given a function $f\in C^\infty(\erren,\erre), f\geq 0,$ we agree to let
\begin{equation}\label{seicinque} \Gamma\ast f =\lim_{m\ttende\infty} \Gamma\ast f_m
\end{equation}
where $f_m=f\varphi_m$ and $\varphi_m\in C_0^\infty (\erren,\erre)$ satisfies
$$\varphi_m=1 \inn D(0,m), \varphi_m=0\inn \erren\meno D(0,m+1)\andd 0\le\varphi_m\le1.$$
Since $0\le\varphi_m\le \varphi_{m+1}$, the sequence $(f_m)_{m\in\enne}$ is nonnegative and increasing, so that 
\eqref{seicinque} is meaningful by Beppo Levi Theorem. It is also easy to recognize that the left hand side of 
\eqref{seicinque} is independent of the choice of the sequence $(\varphi_m)_{m\in\enne}.$ 

The proof of  Theorem \ref{sharp} relies on the following Lemma \ref{lemmaseitre} and Lemma \ref{lemmaseiquattro}. 
\begin{lemma}\label{lemmaseitre} Let $u\in C^\infty(\erren,\erre)$ be such that 
$$u\le 0 \andd \elle u\geq 0\inn \erren.$$
Then
$$u=-\Gamma\ast \elle u + \hat w \ae \inn\erren,$$
where $\hat w\in C^\infty(\erren,\erre),\   \elle \hat w = 0$ and $\hat w \le 0$ in $\erren.$
\end{lemma}
\proof Let $(\varphi_m)_{m\in\enne}$ be a sequence as above and let 
$$f_m=(\elle u)\varphi_m.$$ Then $(f_m)_{m\in\enne}$ is an increasing sequence of nonnegative 
$C_0^\infty(\erren,\erre)$-functions such that 
$$f_m \nearrow f:=\elle u.$$
Define 
\begin{equation}\label{seisette} w_m:= u+\Gamma\ast f_m.
\end{equation}
Then $w_m \in C^\infty(\erren,\erre)$ and
$$\elle(w_m)= \elle u - f_m = f(1 - \varphi_m).$$
Hence,
\begin{equation}\label{seiotto} \elle (w_m)\geq 0 \inn \erren \andd \elle(w_m)=0\inn D(0,m).
\end{equation}
Moreover, since $u\le 0,$
$$\limsup_{x\ttende\infty} w_m(x) \le \limsup_{x\ttende\infty} \Gamma\ast f_m(x)=0.$$
The Maximum Principle of Proposition \ref{maxprinc61} gives
$$w_m\le 0\inn\erren.$$
On the other hand, $(w_m)$ is increasing so that 
$$w:=\lim_{m\ttende\infty} w_m$$ is well defined and satisfies 
\begin{equation}\label{seinove} w_1\le w\le 0.
\end{equation}
This implies $w\in \elleunoloc(\erren)$ and, keeping in mind the second statement in \eqref{seiotto},
$$\int_\erren w \elle^* \varphi\ dx = \lim_{m\ttende \infty} \int_\erren w_m \elle^*\varphi\ dx =0$$
for every $\varphi\in C_0^\infty(\erren).$ As a consequence, since $\elle$ is hypoelliptic, there exists 
$\hat w \in C^\infty(\erren,\erre)$ such that 
$$\hat w = w \ae \andd \elle \hat w=0\inn\erren.$$
Obviously, $\hat w$ also satisfies
$$w_1\le \hat w \le 0\inn\erren.$$
Letting $m$ go to infinity in \eqref{seisette} we obtain 
\begin{equation}\label{seidieci} w = u + \Gamma \ast f,
\end{equation} so that $$u=-\Gamma \ast \elle u + \hat w \ae.$$
\endproof \noindent {\bf{Note.}} \eqref{seidieci} and \eqref{seinove} imply $\Gamma\ast f(x) <\infty$ for every $x\in\erren.$
We complete our prerequisites by proving next lemma.
\begin{lemma}\label{lemmaseiquattro} Let $f\in \cinfty, f\geq 0$ and such that 
\begin{equation}\label{seiundici}\Gamma\ast f\in L^p(\erren) \quad \mbox{ for a suitable } p\in \left[1,1 + \frac{2}{Q-2}\right].
\end{equation}
Then $f\equiv 0.$

Moreover, for every $f\in C_0^\infty(\erren,\erre)$,
$$\Gamma\ast f\in \ellep\quad \mbox{ for every }  p\in \left]1 +\frac{2}{Q-2}, \infty\right[.$$
\end{lemma}
\proof For every $x=(x_1,\ldots,x_n)\in\erren$ define
$$\|x\|= \sum_{j=1}^n |x_j|^{\frac{1}{\sigma_j}}$$
where the $\sigma_j$'s are the exponents related to the dilation $\delta_\lambda$ in \eqref{dilations}. 
Then $x\ttende \|x\|$ is $\delta_\lambda$-homogeneous of degree one:
$$\delta_\lambda(x)\|=\lambda\|x\|\qquad \forall x\in\erren,\forall \lambda>0.$$
Let $\Sigma:=\{x\in\erren\ |\ \|x\|=1\}.$ Since $\Gamma\geq 0$ and, obviously, $\Gamma\nequiv 0$, there exists a (relatively) open subset $\Sigma_0$ of $\Sigma$ such that 
$$\Gamma(x)\geq 2\sigma >0\qquad \forall x\in\Sigma_0,$$
for a suitable $\sigma>0.$ Then, there exists $\rho>0$ such that 
\begin{equation}\label{seitredici}\Gamma(y^{-1}\circ x)\geq \sigma\qquad \forall x\in\Sigma_0,\ \forall y\in\erren\ :\ \|y\|<\rho.
\end{equation}
Consider the {\it open $\delta_\lambda$-cone}
$$K:=\{\delta_\lambda(x)\ |\ x\in \Sigma_0, \  \lambda>0\}.$$
Now, assume by contradiction $f\nequiv 0$. Then there exist a bounded open set $B\neq\emptyset$ such that $f(x)\geq\eps$ for every $x\in B$ and a suitable $\eps>0$. As a consequence, for every $x\in\erren,$ 
\begin{equation}\label{seiquattordici} \begin{split}\Gamma\ast f(x) & \geq
\int_B \Gamma(y^{-1}\circ x) f(y)\ dy \geq \eps \int_B \Gamma(y^{-1}\circ x)\ dy\\&=\eps \|x\|^{2-Q} \int_B \Gamma \left( \left(\delta_{\frac{1}{\|x\|}} (y) \right)^{-1} \circ \delta_{\frac{1}{\|x\|}} (x)\right)\ dy.\end{split}
\end{equation} 
On the other hand, for a suitable $M>1$, $\|\delta_{\frac{1}{\|x\|}} (y)^{-1}\|= \frac{1}{\|x\|} \|y^{-1}\|<\rho$ for every $y\in B$ and $\|x\|\geq M.$ Moreover $\delta_{\frac{1}{\|x\|}} (x) \in \Sigma_0$ if $x\in K.$ Then, by \eqref{seitredici},
$$\Gamma \left( \left( \delta_{\frac{1}{\|x\|}} y \right)^{-1} \circ \delta_{\frac{1}{\|x\|}} (x)\right) \geq \sigma \qquad\forall x\in K,\  \|x\|\geq M\andd \forall y\in B.$$ 
Using this estimate in \eqref{seiquattordici} we get
$$\Gamma\ast f(x)\geq \eps \sigma \|x\|^{2-Q}\qquad \forall x \in K,\ \|x\|\geq M.$$
Therefore:
\begin{equation*}\begin{split} \int_\erren (\Gamma\ast f(x))^p & \geq (\eps \sigma)^p \int_{K\inte \{\|x\|\geq M\}} \|x\|^{p(2-Q)} \ dx  \\&=(\eps \sigma)^p \sum_{k=1}^\infty \int_{K\inte \{M^k \le \|x\| < M^{k+1}\}} \|x\|^{p(2-Q)} \ dx  \\&=\mbox{(using the change of variable $x=\delta_{M^k} (y)$)}
\\&\phantom{=}  (\eps \sigma)^p \left( \int_{K\inte \{ 1\le \|y\|\le M\}} \|y\|^{p(2-Q)}\ dy\right)\sum_{k=1}^\infty M^{k(p(2-Q)+Q)} 
\\&=\infty \mbox{\qquad if \quad }  p\le \frac{Q}{Q-2}= 1+ \frac{2}{Q-2}.\end{split}\end{equation*}
This contradicts the assumption \eqref{seiundici} and proves the first part of the lemma.

To prove the second part we argue as follows. If $f\in C_0^\infty(\erren,\erre)$ then $f\in L^q(\erren)$ for every 
$q\in ]1, \frac{Q}{2}[.$ As a consequence, since $\Gamma$ is $\delta_\lambda$-homogeneous of degree $2-Q$, hence 
$\Gamma\in L^r_{\mathrm{deb}}$ with $r=\frac{Q-2}{Q},$ one has
$$\Gamma\ast f\in L^p(\erren) \mbox{\quad with\quad} \frac{1}{p} =  \frac{1}{r} +  \frac{1}{q} -1 = \frac{1}{q} - \frac{2}{Q}.$$
Since we can choose any $q\in] 1,  \frac{Q}{2}[$, this gives 
 $$\Gamma\ast f \in L^p(\erren) \qquad \forall\  p \in\left] 1 + \frac{2}{Q-2}, \infty\right[.$$
\endproof 
We are ready to prove Theorem \ref{sharp}. Let $u\in \ellep$, with $p\in[1,1+\frac{2}{Q-2}],$ be such that 
$\elle u\geq 0$ in $\erren$ in the weak sense of distributions. We have to prove that $u=0$ a.e. in $\erren.$

We will prove the theorem on the extra assumption $u\in C^\infty(\erren,\erre).$ This restriction can be removed with an approximation  argument like the one used in the proof on Theorem \ref{terzo}. By Theorem \ref{terzo} we already know that $u\le 0$ so that, from Lemma \ref{lemmaseitre}, we get 
$$u=-\Gamma \ast \elle u + \hat w \ae \inn\erren,$$
where $\hat w \in C^\infty(\erren,\erre),\ \elle\hat w =0$ and $\hat w \le 0.$ Then, since $\elle u \geq 0,$ 
$u\le \hat w\le 0.$ Hence, being $u\in\ellep$ with $1\le p \le 1 + \frac{2}{Q-2},$ we also have 
$$\hat w \in \ellep \mbox{\quad for a suitable\ } p\in \left[1, 1 + \frac{2}{Q-2}\right].$$
Theorem \ref{primo} implies $\hat w \equiv 0$, so that 
$$u=-\Gamma\ast \elle u,\inn\erren.$$
From Lemma \ref{lemmaseiquattro} it follows $\elle u=0$, hence $u=0$ in $\erren.$ This completes the proof of the first part of Theorem \ref{sharp}. The second part directly follows from the second part of Lemma \ref{lemmaseiquattro}. 

\section{Proof of Theorem \ref{quinto}}
Let $u\in C^2(\erren,\erre)$ be a classical solution to the equation \eqref{unocinquen} and define 
$$v=F(u)$$ with $F$ given by \eqref{unoquattron}.  Then $v \in C^2(\erren,\erre)$ and, by Lemma \ref{elementarylemma},
\begin{equation*}\begin{split} \elle v &= F'(u)\elle u + F''(u) \gradientea ^2  \\ &= (f(u))^2 + F''(u)\gradientea^2\geq 0.
\end{split}
\end{equation*}
Since $v\in L^1(\erren),$ by Proposition \ref{0}, it follows $\elle v=0,$ i.e.,
$$(f(u))^2 + F''(u)\gradientea^2=0\inn\erren.$$
Being $F''=f'\geq 0,$ from this identity we obtain $f(u)=0,$ hence $u=0$ in $\erren.$

\section{Some examples} \label{examples}
\subsection{} Let $\gi=(\erren,\circ,\delta_\lambda)$ be a stratified Lie group\footnote{We refer to the monograph \cite{BLU} for notions and results recalled in this section.} and let $X_1,\ldots,X_p$ be a basis of the first layer of its Lie algebra. The {\it sub-Laplacian} 
\begin{equation}\label{sublaplacian}
\elle=\sum_{j=1}^p X_j^2 \end{equation}
is left translation invariant on $\gi$ and $\delta_\lambda$-homogeneous of degree two. Then  Theorems \ref{primo}--\ref{quinto} and Corollaries \ref{sesto} and \ref{settimo} apply to $\elle$ in \eqref{sublaplacian}.
\subsection{} Let $\gi$ as above and consider in $\erre^{n+1}:=\erre_x^n\times\erre_t$ the {\it heat-type operator} 
\begin{equation}\label{calore} \acca:= \sum_{j=1}^p X_j^2 -\partial_t.
\end{equation} 
This operator is left translation invariant and homogeneous of degree two with respect to the stratified 
Lie group
$$\hat\gi=\gi\oplus \erre=(\erre^{n+1},\hat\circ,\hat\delta_\lambda),$$
where $\hat\circ$ and $\hat\delta_\lambda$ are defined as follow
\begin{equation*}\begin{split}&(x,t) \hat\circ(x',t') = (x\circ x', t+t')\\&\hat\delta_\lambda(x,t)=(\delta_\lambda(x), \lambda^2 t).
\end{split}
\end{equation*} 
The {\it homogeneous dimension} of $\hat\gi$ is 
$$\hat Q = Q+2$$
being $Q$ the homogeneous dimension of $\gi.$\\
To the operator $\acca$ in \eqref{calore}  Theorems \ref{primo}--\ref{quinto} and Corollaries \ref{sesto} and \ref{settimo} apply.
\subsection{}
Let us consider in $\erre^{n+1}:=\erre^n_x\times\erre_t$ the {\it Kolmogorov-type operators} 
\begin{equation}\label{kolmogorov} \elle = \dive (A\nabla) + \langle Bx,\nabla \rangle - \partial_t,\end{equation}
where $A$ and $B$ are constant $n\times n$ real matrices, $A$ symmetric and $\geq 0$.

Define
$$E(s):=\exp(-sB),\quad s\in\erre.$$
Then the operator $\elle$ in \eqref{kolmogorov} is left translation invariant on the Lie group
$$\mathbb{K}=(\erre^{n+1},\circ)$$
with composition law
\begin{equation}\label{compositionlaw}(x,t)\circ (x',t')=(x'+E(t')x, t+t').\end{equation}
The Lebesque measure is both left and right invariant on $\mathbb{K}$ if and only if
\begin{equation}\label{traccia}\mathrm{trace} (B)=0.\end{equation}
Moreover, if we assume 
\begin{equation}\label{ipo}C(t):=\int_0^t E(s) A E^T(s)\ ds >0\qquad\forall\ t>0,\end{equation}
then $\elle$ is hypoelliptic (see e.g. \cite{lanconelli_polidoro_1994}, see also \cite[Sections 4.1.3, 4.3.4]{BLU}).
Then, under the assumptions \eqref{traccia} and \eqref{ipo}, Theorems \ref{primo}, \ref{secondo}, \ref{terzo}, \ref{quinto} and Corollaries \ref{sesto} and \ref{settimo}  apply to $\elle$ in \eqref{kolmogorov}.
On the other hand, if the matrix $B$ takes the particular block form fixed in \cite{lanconelli_polidoro_1994}, then there exists a family of dilations $(\delta_\lambda)_{\lambda > 0}$ in $\erre^{n+1}$ making
$$\mathbb{K}=(\erre^{n+1}, \circ, \delta_\lambda)$$
a homogeneous Lie group and the operator $\elle$ in \eqref{kolmogorov} is $\delta_\lambda$-homogeneous of degree two. Therefore under this extra assumption, also 
Theorem \ref{sharp} apply to  $\cal L$.

\begin{remark}\label{priola} {\rm  Consider the {\it stationary} counterpart of $\elle$ in \eqref{kolmogorov}, i.e., the degenerate Ornstein-Uhlenbeck operator 
\begin{equation}\label{ou} \elle = \dive (A\nabla) + \langle Bx,\nabla \rangle.\end{equation}
Priola and Zabczyk in \cite[Theorem 3.1]{priola_zabczyk_2004} proved that  $\elle_0$ has the {\it $L^\infty$-Liouville property if and only if}}
\begin{equation*} \mathrm{Re} (\lambda) \le 0 \mbox{ \it for every $\lambda$ eigenvalue of $B$.}
\end{equation*} 
{\rm Then, if $B$ has an eigenvalue with real part strictly positive, there exists a bounded solution $v$ to $\elle_0 v=0$ in $\erren$ which is not constant. 
Hence 
$$u(x,t)=v(x)$$ is a bounded nonconstant solution to
$$\elle u=0 \inn\erre^{n+1}.$$
Thus {\it $\elle$ does not have the $L^\infty$-Liouville property.} An explicit example is given by the operator $\elle$
in Remark \ref{noinfty}, which can be written as in \eqref{kolmogorov} by taking 

\begin{displaymath}
A=\left(\begin{array}{cc}1 & 0 \\ 0& 0\end{array}
\right) \andd B=\left(\begin{array}{cc}1 & -\dfrac{1}{2} \\ -\dfrac{1}{2} & -1\end{array}
\right).\end{displaymath}
The eigenvalues of $B$ are $ -\frac{\sqrt{3}}{2}$ and  $ \frac{\sqrt{3}}{2}$. Moreover 
$$\mathrm{trace}(B)=0,$$
so that the Lebesgue measure is both left and right invariant w.r.t. the composition law in \eqref{compositionlaw}.
Finally, \eqref{ipo} can be 
verified by a direct computation or simply recognizing that, letting 
$$X=\partial_{x_1}\andd Y=\left(x_1 - \dfrac{1}{2}x_2\right)\partial_{x_1}+\left(\dfrac{1}{2}x_1 - x_2\right)\partial_{x_2}-\partial_t,$$
the hypoellipticity H\"ormander rank condition $$\mathrm{Lie} \{X,Y\}(x,t)=3\qquad\forall\ (x,t)\in\erre^3$$
is satisfied.
Then: {\it the operator $\elle$ in Remark \ref{noinfty} has the $L^p$-Liouville property for every $p\in[0,\infty[$, but it has not the $L^\infty$-Liouville property.}
}\end{remark} 
\subsection{} 
Let us consider the operator in $\erre^n$, $n\geq 3$,
\begin{equation}\label{cheneso}\elle=\partial_{x_1}^2 + u_1(x_1)\partial_{x_2}+\cdots+u_{n-1}(x_1)\partial_{x_n}\end{equation} where $\{u_1,\ldots,u_{n-1}\}$ is a real fundamental system of solutions of the ODE $P(u)=0$, being 
$$P(u):=u^{(n-1)}+ a_{n-2}u^{(n-2)}+\cdots+a_1u^{(1)} +a_0u,$$
with $a_{0}, \ldots, a_{n-2}\in\erre$ and $ a_{n-2}=1.$ In \cite{bonfiglioli_lanconelli_2012} it is proved that
$\elle$ is hypoelliptic and left translation invariant  on a Lie group $\gi(P)=(\erren,\circ)$ which, in \cite{bonfiglioli_lanconelli_2012}, is called {\it $P$-group.} Due to the condition  $a_{n-2}=1$, the Lebesgue measure 
is both left and right invariant on $\gi(P)$. Then Theorems \ref{primo}, \ref{secondo}, \ref{terzo}, \ref{quinto} and 
Corollaries \ref{sesto} and \ref{settimo} apply to the operator $\elle$ in \eqref{cheneso}.
\begin{remark}{\rm An explicit example of an operator $\elle$ as in \eqref{cheneso} is given by the {\it Mumford operator}
$$\mathcal{M}:=\partial_{x_1}^2 + \cos x_1 \partial_{x_2} +\sen x_1 \partial_{x_3} \inn \erre^3$$
which is left invariant on $\gi(P)$ with $P$ given by 
$$P(u)=u''+u.$$
}\end{remark}
\begin{remark}{\rm Theorems \ref{primo}, \ref{secondo}, \ref{terzo}, \ref{quinto} and 
Corollaries \ref{sesto} and \ref{settimo}  also apply to the {\it evolution counterpart} of the operator $\elle$ in \eqref{cheneso}, i.e., to 
\begin{equation*}\elle - \partial_t =\partial_{x_1}^2 + u_1(x_1)\partial_{x_2}+\cdots+u_{n-1}(x_1)\partial_{x_n} - \partial_t .\end{equation*} This operator is hypoelliptic and left translation invariant w.r.t.  the group composition law 
$$(x,t)\hat\circ(x',t')=(x\circ x', t+t'),$$ where 
$\circ$ is the composition law of $\gi(P)$. More formally, $\elle - \partial_t$ is left invariant on $$\gi(P)\oplus\erre.$$

}\end{remark} 

\section*{Acknowledgement}
A preliminary version of this paper
was presented by the first author
in the  {\it Bruno Pini Analysis Mathematical Seminars}
of the Department of Mathematics of Bologna (\cite{kogoj_BPMAS}).
\bibliographystyle{alpha}
\bibliography{bibliografia}
\end{document}